\documentclass[a4paper,12pt]{article}

 \usepackage{latexsym,amssymb,amsmath}

 \newtheorem{theorem}{Theorem}
 \newtheorem{lemma}{Lemma}
 \newtheorem{proposition}{Proposition}

\begin{document}

\title{Fueter's theorem: the saga continues\thanks{accepted for publication in \textit{Journal of Mathematical Analysis and Applications}}}

\author{Dixan Pe\~na Pe\~na$^{\star,1}$ and Frank Sommen$^{\star\star,2}$}

\date{\normalsize{$^\star$Department of Mathematics, Aveiro University,\\3810-193 Aveiro, Portugal\\
$^{\star\star}$Department of Mathematical Analysis, Ghent University,\\9000 Gent, Belgium}\\\vspace{0.4cm}
\small{$^1$e-mail: dixanpena@ua.pt; dixanpena@gmail.com\\
$^2$e-mail: fs@cage.ugent.be}} 

\maketitle

\begin{abstract}
\noindent In this paper is extended the original theorem by Fueter-Sce (assigning an $\mathbb R_{0,m}$-valued monogenic function to a $\mathbb C$-valued holomorphic function) to the higher order case. We use this result to prove Fueter's theorem with an extra monogenic factor $P_k(x_0,\underline x)$.\vspace{0.2cm}\\
\textit{Keywords}: Fueter's theorem, CK-extension, Almansi-Fischer decomposition.\vspace{0.1cm}\\
\textit{Mathematics Subject Classification}: 30G35.
\end{abstract}

\section{Introduction and preliminaries}

Let $\mathbb{R}_{0,m}$ be the $2^m$-dimensional real Clifford algebra constructed over the orthonormal basis $(e_1,\ldots,e_m)$ of the Euclidean space $\mathbb R^m$ (see \cite{Cl}). The multiplication in $\mathbb{R}_{0,m}$ is determined by the relations 
\[e_je_k+e_ke_j=-2\delta_{jk},\quad j,k=1,\dots,m,\]
where $\delta_{jk}$ is the Kronecker delta. A basis for the algebra $\mathbb{R}_{0,m}$ is then given by the elements 
\[e_A=e_{j_1}\cdots e_{j_k},\]
where $A=\{j_1,\dots,j_k\}\subset\{1,\dots,m\}$ is such that $j_1<\dots<j_k$. For the empty set $\emptyset$, we put $e_{\emptyset}=1$, the latter being the identity element.

Any Clifford number $a\in\mathbb{R}_{0,m}$ may thus be written as 
\[a=\sum_Aa_Ae_A,\quad a_A\in\mathbb R,\]
and its conjugate $\overline a$ is given by
\[\overline a=\sum_Aa_A\overline e_A,\quad\overline e_A=\overline e_{j_k}\dots\overline e_{j_1},\;\overline e_j=-e_j,\;j=1,\dots,m.\]
Observe that $\mathbb R^{m+1}$ may be naturally embedded in $\mathbb R_{0,m}$ by associating to any element $(x_0,x_1,\ldots,x_m)\in\mathbb R^{m+1}$ the paravector $x=x_0+\underline x=x_0+\sum_{j=1}^mx_je_j$.

Let us recall that an $\mathbb R_{0,m}$-valued function $f$ defined and continuously differentiable in an open set $\Omega$ of $\mathbb R^{m+1}$, is said to be (left) monogenic in $\Omega$ if $\partial_xf=0$ in $\Omega$, where
\[\partial_x=\partial_{x_0}+\partial_{\underline x}\]
is the generalized Cauchy-Riemann operator in $\mathbb R^{m+1}$ and 
\[\partial_{\underline x}=\sum_{j=1}^me_j\partial_{x_j}\]
is the Dirac operator in $\mathbb R^m$. Null-solutions of $\partial_{\underline x}$ are also called monogenic functions. 

These functions are a fundamental object in Clifford analysis and may be considered as a natural generalization to higher dimension of the holomorphic functions in the complex plane (see e.g. \cite{BDS}). 

It is worth remarking that $\partial_{\underline x}$ and $\partial_x$ factorize the Laplacian, i.e. 
\begin{align*}
\Delta_{\underline x}&=\sum_{j=1}^m\partial_{x_j}^2=-\partial_{\underline x}^2,\\
\Delta_x&=\partial_{x_0}^2+\Delta_{\underline x}=\partial_x\overline\partial_x=\overline\partial_x\partial_x,
\end{align*}
and therefore every monogenic function is also harmonic.

An interesting fact about the monogenic functions was first observed by Fueter (see \cite{F}) in the setting of quaternionic analysis: it is possible to generate monogenic functions using holomorphic functions in the upper half of the complex plane. This fact, which is now known as Fueter's theorem, was later extended to the case of $\mathbb R_{0,m}$-valued functions by Sce \cite{Sce} ($m$ odd) and Qian \cite{Q} ($m$ even). For further works on this topic we refer the reader to \cite{KQS,LaLe,LaRa,D,DS,DS2,DS3,QS,S}.

Let $f(z)=u(t_1,t_2)+iv(t_1,t_2)$ ($z=t_1+it_2$) be a holomorphic function in some open subset $\Xi\subset\mathbb C^+=\{z\in\mathbb C:\;t_2>0\}$; and let $P_k(\underline x)$ be a homogeneous monogenic polynomial of degree $k$ in $\mathbb R^m$, i.e.
\begin{alignat*}{2}
\partial_{\underline x}P_k(\underline x)&=0,&\quad\underline x&\in\mathbb R^m,\\
P_k(t\underline x)&=t^kP_k(\underline x),&\quad t&\in\mathbb R.
\end{alignat*}
Put $\underline\omega=\underline x/r$, with $r=\vert\underline x\vert=\sqrt{-\underline x^2}$. In this paper we will focus on the following generalization of Fueter's theorem obtained in \cite{S}.

\begin{theorem}\label{Som-Fue-thm}
If $m$ is odd, then the function
\begin{equation*}
\Delta_x^{k+\frac{m-1}{2}}\bigl[\bigl(u(x_0,r)+\underline\omega\,v(x_0,r)\bigr)P_k(\underline x)\bigr]
\end{equation*}
is monogenic in $\widetilde\Omega=\{x\in\mathbb R^{m+1}:\;(x_0,r)\in\Xi\}$.
\end{theorem}
\noindent
\textbf{Remark:} The case $k=0$ corresponds to Sce's result. 

The purpose of this paper is to prove that Theorem \ref{Som-Fue-thm} is still valid if we replace $P_k(\underline x)$ by a homogeneous monogenic polynomial $P_k(x_0,\underline x)$ of degree $k$ in $\mathbb R^{m+1}$ (this problem arose from discussions between Qian and Sommen). In order to attain this goal, we will first prove an extension of Fueter's theorem that uses $\mathbb C$-valued functions satisfying the equation 
\[\partial_{\bar z}\Delta_z^pf(z)=0,\quad p\in\mathbb N_0,\]
as initial functions instead of the usual holomorphic functions. Here and throughout $\partial_{\bar z}$ and $\Delta_z$ denote, respectively, the well-known Cauchy-Riemann operator 
\[\partial_{\bar z}=\frac{1}{2}(\partial_{t_1}+i\partial_{t_2})\]
and the Laplacian in two dimensions
\[\Delta_z=\partial_{t_1}^2+\partial_{t_2}^2.\]

\section{A higher order version of the original\\Fueter's theorem}

We begin this section with two essential lemmas. The proof of Lemma \ref{lemma1} may be found in \cite{D,DS}. 

\begin{lemma}\label{lemma1}
Suppose that $f(t_1,\dots,t_d)$ is an $\mathbb R$-valued infinitely differentiable function in an open set of $\mathbb R^d$ and that $D_{t_j}(n)$ and $D^{t_j}(n)$, $n\in\mathbb N_0$, are differential operators defined by  
\begin{align*}
D_{t_j}(n)\{f\}&=\left(\frac{1}{{t_j}}\,\partial_{t_j}\right)^n\{f\},\quad\quad\quad\;\;\,D_{t_j}(0)\{f\}=f,\\
D^{t_j}(n)\{f\}&=\partial_{t_j}\left(\frac{D^{t_j}(n-1)\{f\}}{{t_j}}\right),\;\,D^{t_j}(0)\{f\}=f,
\end{align*}
$j=1,\dots,d$. Then one has
\begin{itemize}
\item[{\rm(i)}] $\partial_{t_j}^2D_{t_j}(n)\{f\}=D_{t_j}(n)\{\partial_{t_j}^2f\}-2n D_{t_j}(n+1)\{f\}$, 
\item[{\rm(ii)}] $\partial_{t_j}^2D^{t_j}(n)\{f\}=D^{t_j}(n)\{\partial_{t_j}^2f\}-2nD^{t_j}(n+1)\{f\}$,
\item[{\rm(iii)}] $D^{t_j}(n)\{\partial_{t_j}f\}=\partial_{t_j}D_{t_j}(n)\{f\}$,  
\item[{\rm(iv)}] $D_{t_j}(n)\{\partial_{t_j}f\}-\partial_{t_j}D^{t_j}(n)\{f\}=2n/t_j\,D^{t_j}(n)\{f\}$.
\end{itemize}
\end{lemma}

\begin{lemma}\label{mainlemma}
Suppose that $g$ is an $\mathbb R$-valued infinitely differentiable function in an open set of $\mathbb R^2_+=\{(t_1,t_2)\in\mathbb R^2:\;t_2>0\}$.
Then for $n\in\mathbb N_0$, we have that  
\begin{itemize}
\item[{\rm(i)}] $\displaystyle{\Delta_x^n\big(g(x_0,r)P_k(\underline x)\big)=\left(\sum_{j=0}^nd_{k,m}(j)\binom{n}{j}D_r(j)\Delta_z^{n-j}g(x_0,r)\right)P_k(\underline x)}$,
\item[{\rm(ii)}] $\displaystyle{\Delta_x^n\big(g(x_0,r)\,\underline\omega\,P_k(\underline x)\big)=\left(\sum_{j=0}^nd_{k,m}(j)\binom{n}{j}D^r(j)\Delta_z^{n-j}g(x_0,r)\right)\underline\omega\,P_k(\underline x)}$,
\end{itemize}
where 
\begin{align*}
d_{k,m}(j)&=(2k+m-1)(2k+m-3)\cdots(2k+m-(2j-1)),\\d_{k,m}(0)&=1.
\end{align*}
\end{lemma}
\textit{Proof.} It is easily seen that 
\begin{align*}
\partial_{\underline x}g&=\sum_{j=1}^me_j\partial_{x_j}g=\sum_{j=1}^me_j(\partial_rg)(\partial_{x_j}r)\\
&=\sum_{j=1}^me_j(\partial_rg)\frac{x_j}{r}=\underline\omega\partial_rg,
\end{align*}
\[\Delta_{\underline x}\,\underline\omega=-\partial_{\underline x}^2\,\underline\omega=(m-1)\partial_{\underline x}\left(\frac{1}{r}\right)=-\frac{(m-1)}{r^2}\,\underline\omega,\]
\begin{align*}
\Delta_xg&=\partial_{x_0}^2g+\Delta_{\underline x}g=\partial_{x_0}^2g-\partial_{\underline x}(\underline\omega\partial_rg)\\
&=\partial_{x_0}^2g+\partial_r^2g+\frac{m-1}{r}\,\partial_rg.
\end{align*}
Therefore  
\[\hspace{-1.4cm}\Delta_x(gP_k)=(\Delta_xg)P_k+2\sum_{j=1}^m(\partial_{x_j}g)(\partial_{x_j}P_k)+g(\Delta_{\underline x}P_k)\]
\[\hspace{1.8cm}=\left(\partial_{x_0}^2g+\partial_r^2g+\frac{m-1}{r}\,\partial_rg\right)P_k+2\frac{\partial_rg}{r}\sum_{j=1}^mx_j\partial_{x_j}P_k\]
\[\hspace{-0.7cm}=\left(\partial_{x_0}^2g+\partial_r^2g+\frac{2k+m-1}{r}\,\partial_rg\right)P_k\]
\begin{equation}\label{eq1}
\hspace{0.3cm}=\left(\partial_{x_0}^2g+\partial_r^2g+(2k+m-1)D_r(1)\{g\}\right)P_k
\end{equation}
and 
\[\hspace{-2.9cm}\Delta_x(g\underline\omega P_k)=(\Delta_{\underline x}\,\underline\omega)gP_k+2\sum_{j=1}^m(\partial_{x_j}\underline\omega)(\partial_{x_j}(gP_k))+\underline\omega\Delta_x(gP_k)\]
\[\hspace{1.55cm}=-\frac{(m-1)}{r^2}\,g\underline\omega P_k+2\sum_{j=1}^m\left(\frac{e_j}{r}-\frac{x_j}{r^2}\,\underline\omega\right)\left(\frac{x_j}{r}\,(\partial_rg)P_k+g(\partial_{x_j}P_k)\right)\]
\[\hspace{-1.5cm}+\left(\partial_{x_0}^2g+\partial_r^2g+\frac{2k+m-1}{r}\,\partial_rg\right)\underline\omega P_k\]
\[\hspace{-0.99cm}=\left(\partial_{x_0}^2g+\partial_r^2g+(2k+m-1)\left(\frac{\partial_rg}{r}-\frac{g}{r^2}\right)\right)\underline\omega P_k\]
\begin{equation}\label{eq2}
\hspace{-1.75cm}=\left(\partial_{x_0}^2g+\partial_r^2g+(2k+m-1)D^r(1)\{g\}\right)\underline\omega P_k,
\end{equation}
where we have also used Euler's theorem for homogeneous functions.

The proof now follows by induction using equalities (\ref{eq1})-(\ref{eq2}) together with statements (i)-(ii) of Lemma \ref{lemma1}.\hfill$\square$\vspace{0.22cm}

The previous lemma allows us to obtain a method to generate functions on $\mathbb R^{m+1}$ that satisfy the equation 
\begin{equation}\label{pseudopharm2}
\Delta_x^pF(x_0,\underline x)=0,
\end{equation} 
using functions fulfilling an equation of the same type in $\mathbb R^2$ 
\begin{equation}\label{pseudopharm1}
\Delta_z^pg(t_1,t_2)=0.
\end{equation}

\begin{proposition}
Let $m$ be odd. If $g$ is an $\mathbb R$-valued solution of (\ref{pseudopharm1}) in $\Xi\subset\mathbb R^2_+$, then
\[\Delta_x^{k+\frac{m-1}{2}}\big(g(x_0,r)P_k(\underline x)\big)\]
and
\[\Delta_x^{k+\frac{m-1}{2}}\big(g(x_0,r)\,\underline\omega\,P_k(\underline x)\big)\]
satisfy the equation (\ref{pseudopharm2}) in $\widetilde\Omega=\{x\in\mathbb R^{m+1}:\;(x_0,r)\in\Xi\}$.
\end{proposition}
\textit{Proof.} We first observe that if $m$ is odd, then each factor which appears in the definition of $d_{k,m}(j)$ is even. Now, by Lemma \ref{mainlemma}, we have that
\begin{multline*}
\Delta_x^{p+k+\frac{m-1}{2}}\big(g(x_0,r)P_k(\underline x)\big)\\=\left(\sum_{j=0}^{p+k+\frac{m-1}{2}}d_{k,m}(j)\binom{p+k+\frac{m-1}{2}}{j}D_r(j)\Delta_z^{p+k+\frac{m-1}{2}-j}g(x_0,r)\right)P_k(\underline x),
\end{multline*}
\begin{multline*}
\Delta_x^{p+k+\frac{m-1}{2}}\big(g(x_0,r)\,\underline\omega\,P_k(\underline x)\big)\\=\left(\sum_{j=0}^{p+k+\frac{m-1}{2}}d_{k,m}(j)\binom{p+k+\frac{m-1}{2}}{j}D^r(j)\Delta_z^{p+k+\frac{m-1}{2}-j}g(x_0,r)\right)\underline\omega\,P_k(\underline x).
\end{multline*}
Clearly, the first $k+(m+1)/2$ terms in the above equalities vanish since $g$ is by hypothesis a solution of (\ref{pseudopharm1}). Finally, note that $2k+m-(2j-1)\le0$ for $j\ge k+(m+1)/2$ and therefore $d_{k,m}(j)=0$ for $j\ge k+(m+1)/2$.\hfill$\square$

\begin{theorem}\label{msFThm}
Let $f(z)=u(t_1,t_2)+iv(t_1,t_2)$ be a $\mathbb C$-valued function satisfying in some open subset $\Xi\subset\mathbb C^+$ the equation
\[\partial_{\bar z}\Delta_z^pf(z)=0,\quad p\in\mathbb N_0.\]
If $m$ is odd, then the function
\begin{equation*}
\Delta_x^{p+k+\frac{m-1}{2}}\bigl[\bigl(u(x_0,r)+\underline\omega\,v(x_0,r)\bigr)P_k(\underline x)\bigr]
\end{equation*}
is monogenic in $\widetilde\Omega=\{x\in\mathbb R^{m+1}:\;(x_0,r)\in\Xi\}$.
\end{theorem}
\textit{Proof.} By Lemma \ref{mainlemma}, we get that 
\begin{multline*}
\Delta_x^{p+k+\frac{m-1}{2}}\bigl[\bigl(u(x_0,r)+\underline\omega\,v(x_0,r)\bigr)P_k(\underline x)\bigr]\\
=(2k+m-1)!!\binom{p+k+\frac{m-1}{2}}{k+\frac{m-1}{2}}\bigl(A(x_0,r)+\underline\omega\,B(x_0,r)\bigr)P_k(\underline x),
\end{multline*}
with 
\[A=D_r\left(k+\frac{m-1}{2}\right)\{\Delta_z^pu\},\]
\[B=D^r\left(k+\frac{m-1}{2}\right)\{\Delta_z^pv\}.\]
It thus remains to prove that $A$ and $B$ satisfy the Vekua-type system (see \cite{LB,D,S2,S3})
\begin{equation*}
\left\{\begin{array}{ll}\partial_{x_0}A-\partial_rB&=\displaystyle{\frac{2k+m-1}{r}}\,B\\\partial_{x_0}B+\partial_rA&=0.\end{array}\right. 
\end{equation*}
In order to do that, it will be necessary to use the fact that $u$ and $v$ satisfy in $\Xi$ the system
\begin{equation*}
\left\{\begin{array}{ll}\partial_{t_1}\Delta_z^pu-\partial_{t_2}\Delta_z^pv&=0\\\partial_{t_1}\Delta_z^pv+\partial_{t_2}\Delta_z^pu&=0,\end{array}\right. 
\end{equation*}
and statements (iii)-(iv) of Lemma \ref{lemma1}. 

Indeed, 
\[\begin{split}
\partial_{x_0}A-\partial_rB&=D_r\left(k+\frac{m-1}{2}\right)\{\partial_{x_0}\Delta_z^pu\}-\partial_rD^r\left(k+\frac{m-1}{2}\right)\{\Delta_z^pv\}\\
&=D_r\left(k+\frac{m-1}{2}\right)\{\partial_r\Delta_z^pv\}-\partial_rD^r\left(k+\frac{m-1}{2}\right)\{\Delta_z^pv\}\\
&=\frac{2k+m-1}{r}\,D^r\left(k+\frac{m-1}{2}\right)\{\Delta_z^pv\}\\
&=\frac{2k+m-1}{r}\,B
\end{split}\]
and 
\[\begin{split}
\partial_{x_0}B+\partial_rA&=D^r\left(k+\frac{m-1}{2}\right)\{\partial_{x_0}\Delta_z^pv\}+\partial_rD_r\left(k+\frac{m-1}{2}\right)\{\Delta_z^pu\}\\
&=D^r\left(k+\frac{m-1}{2}\right)\{\partial_{x_0}\Delta_z^pv\}+D^r\left(k+\frac{m-1}{2}\right)\{\partial_r\Delta_z^pu\}\\
&=D^r\left(k+\frac{m-1}{2}\right)\{\partial_{x_0}\Delta_z^pv+\partial_r\Delta_z^pu\}\\
&=0,
\end{split}\]
which completes the proof.\hfill$\square$

\section{Fueter's theorem with an extra monogenic factor $P_k(x_0,\underline x)$}

Before starting the proof of the main theorem, we need to recall two basic results of Clifford analysis (see \cite{BDS,DSS,MR}). 

\begin{theorem}[CK-extension theorem]\label{CK}
Every function $g(\underline x)$ analytic in $\mathbb R^m$ has a unique monogenic extension $\mathsf{CK}[g]$ to $\mathbb R^{m+1}$, which is given by
\begin{equation*}
\mathsf{CK}[g(\underline x)](x)=\sum_{j=0}^\infty\frac{(-x_0)^j}{j!}\,\partial_{\underline x}^jg(\underline x).
\end{equation*}
\end{theorem}
Let us denote by $\mathsf{P}(k)$, $k\in\mathbb N$, the set of all $\mathbb R_{0,m}$-valued homogeneous polynomials of degree $k$ in $\mathbb R^m$. It contains the important subspace $\mathsf{M}^+(k)$ consisting of all homogeneous monogenic polynomials of degree $k$.

\begin{theorem}[Almansi-Fischer decomposition]\label{Fischer}
Let $k\in\mathbb N$. Then
\[\mathsf{P}(k)=\bigoplus_{n=0}^k\underline x^n\mathsf{M}^+(k-n).\]
\end{theorem}
We are now ready to prove the final result.

\begin{theorem}
Let $f(z)=u(t_1,t_2)+iv(t_1,t_2)$ be a $\mathbb C$-valued holomorphic function in some open subset $\Xi\subset\mathbb C^+$ and assume that $P_k(x_0,\underline x)$ is a homogeneous monogenic polynomial of degree $k$ in $\mathbb R^{m+1}$.
If $m$ is odd, then the function
\begin{equation*}
\Delta_x^{k+\frac{m-1}{2}}\bigl[\bigl(u(x_0,r)+\underline\omega\,v(x_0,r)\bigr)P_k(x_0,\underline x)\bigr]
\end{equation*}
is monogenic in $\widetilde\Omega=\{x\in\mathbb R^{m+1}:\;(x_0,r)\in\Xi\}$.
\end{theorem}
\textit{Proof.} It is clear from Theorem \ref{CK} that $P_k(x_0,\underline x)=\mathsf{CK}\left[P_k(0,\underline x)\right](x)$. By Theorem \ref{Fischer}, there exist unique $P_{k-n}(\underline x)\in\mathsf{M}^+(k-n)$ such that
\[P_k(x_0,\underline x)=\sum_{n=0}^k\mathsf{CK}\left[\underline x^nP_{k-n}(\underline x)\right](x).\]
Thus, it suffices to prove the monogenicity of 
\[\Delta_x^{k+\frac{m-1}{2}}\bigl[\bigl(u(x_0,r)+\underline\omega\,v(x_0,r)\bigr)\mathsf{CK}\left[\underline x^nP_{k-n}(\underline x)\right](x)\bigr],\quad n=0,\dots,k.\]
As
\[\mathsf{CK}\left[\underline x^nP_{k-n}(\underline x)\right](x)=\sum_{j=0}^k\frac{(-x_0)^j}{j!}\,\partial_{\underline x}^j\big(\underline x^nP_{k-n}(\underline x)\big)\] 
and on account of the equality
\[\partial_{\underline x}\big(\underline x^nP_{k-n}(\underline x)\big)=\left\{\begin{array}{ll}-(2k+m-n-1)\underline x^{n-1}P_{k-n}(\underline x),&\text{if}\;n\;\text{odd},\\-n\underline x^{n-1}P_{k-n}(\underline x),&\text{if}\;n\;\text{even},\end{array}\right.\]
we may conclude that $\mathsf{CK}\left[\underline x^nP_{k-n}(\underline x)\right](x)$ is of the form 
\[\mathsf{CK}\left[\underline x^nP_{k-n}(\underline x)\right](x)=\left(\sum_{j=0}^kc_jx_0^j\underline x^{n-j}\right)P_{k-n}(\underline x),\quad c_j\in\mathbb R.\]
Therefore
\[\mathsf{CK}\left[\underline x^nP_{k-n}(\underline x)\right](x)=\big(U(x_0,r)+\underline\omega\,V(x_0,r)\big)P_{k-n}(\underline x),\]
where $U$ and $V$ are $\mathbb R$-valued homogeneous polynomial of degree $n$ in the variables $x_0$ and $r$. Its corresponding $\mathbb C$-valued function $g(z)=U(t_1,t_2)+iV(t_1,t_2)$ clearly satisfies
\[\partial_{\bar z}^{n+1}g(z)=0,\quad z\in\mathbb C,\] 
whence
\[\partial_{\bar z}^{n+1}\big(f(z)g(z)\big)=0,\quad z\in\Xi,\]
i.e. $f(z)g(z)$ is $(n+1)$-holomorphic in $\Xi$. It then follows that  
\[\partial_{\bar z}\Delta_z^n\big(f(z)g(z)\big)=0,\quad z\in\Xi.\]
The proof now follows by using Theorem \ref{msFThm}.\hfill$\square$

\subsection*{Acknowledgment}

D. Pe\~na Pe\~na was supported by a Post-Doctoral Grant of \emph{Funda\c{c}\~ao para a Ci\^encia e a Tecnologia}, Portugal (grant number: SFRH/BPD/45260/2008).

\end{document}